\documentclass[10pt,leqno]{article}

\usepackage{amsmath}
\usepackage{latexsym}
\usepackage{graphicx}
\makeatletter
\def\@maketitle{\newpage
    \null
    \vskip .8truein
    \begin{center}%
     {\bf \@title \par}%
     \vskip 1.5em
     {\small
      \lineskip .5em
      \begin{tabular}[t]{c}\@author
      \end{tabular}\par}%
    \end{center}%
    \par
    \vskip .4truein}
\@addtoreset{equation}{section} \@addtoreset{theorem}{section}
\@addtoreset{lemma}{section} \@addtoreset{proposition}{section}
\@addtoreset{definition}{section} \@addtoreset{corollary}{section}
\@addtoreset{remark}{section}

\let\s=\sigma

\let\ol=\overline
\newcommand{\re}{{{I\!\!R}}}

\def\R{{\bf R}}
\let\a=\alpha
\let\d=\delta
\let\de=\partial

\let\l=\lambda
\let\vf=\varphi

\newtheorem{theorem}{Theorem}[section]
\newtheorem{lemma}{Lemma}[section]
\newtheorem{proposition}{Proposition}[section]

\DeclareMathOperator{\argmax}{argmax}

\def\proof{\list{}{\setlength{\leftmargin}{0pt}
                      \parskip=0pt\parsep=0pt\listparindent=2em
                      \itemindent=0pt}\item[]\futurelet\testchar\@maybe}

\def\@maybe{\ifx[\testchar \let\next\@Opt
          \else \let\next\@NoOpt \fi \next}
\def\@Opt[#1]{{\it Proof of #1.\ }}\def\@NoOpt{{\it Proof.\ }}
\begin{document}

\title{%\Large \bf 
Comparison principles for subelliptic equations of Monge-Amp\`ere type}

%%%%%%%%%%%%%%%%%%%%%%%%%%%%%%%%%%%%%%%%%%%%%%%%%%%
%\footnote{\date{\today}}\\

\author{%{\large \sc 
Martino Bardi, 
Paola Mannucci\\
%\thanks{
%\address{
Dipartimento di Matematica P. e A.,
Universit\`a di Padova\\
 via Trieste 63, 35121 Padova, Italy}

\maketitle

%\vskip-1.in

\begin{abstract}
We present two comparison principles for viscosity sub- and supersolutions 
of Monge-Amp\`ere-type equations associated to a family of vector fields.
In particular, we obtain the uniqueness of a viscosity solution to the Dirichlet problem for the equation of prescribed horizontal Gauss curvature in a Carnot group. 
%\footnote{\date{\today}}\\
\end{abstract}

%\vskip 1cm

%	\noindent {\bf Keywords}:  viscosity solution, comparison principle, Monge-Amp\`ere % type 
%	equation,
%	subelliptic equation, Carnot groups, Gauss curvature.
%\noindent  {\bf AMS Subject classification:} 35J70, 35H20,
%35J60.
%\footnote{\date{\today}}\\

%\date{\today}

%\newpage
%%%%%%%%%%%%%%%%%%%%%%%%%%%%%%%%%%%%%
%Introduzione 
%%%%%%%%%%%%%%%%%%%%%%%%%%%%%%%%%%%%%%
\section*{Introduction}
We consider fully nonlinear partial differential equations of the form
\begin{equation}
\label{MAh}
- \det (D_ {\mathcal X}^2 u)+ H(x, u, D_{\mathcal X}u)=0, \ \mbox{in}\ \Omega,
\end{equation}
where $\Omega\subseteq\re^n$ is open 
and bounded, $D_{\mathcal X}u$ denotes the gradient of $u$ with respect to  a given family of $C^{1,1}$ vector fields $X_1, ..., X_m$,
%\[
$
D_{\mathcal X}u := (X_1u, ..., X_mu),
$
%\]
$D^2_{\mathcal X}u$ denotes the symmetrized Hessian matrix of $u$ with respect to  the same vector fields
\[
(D^2_{\mathcal X}u)_{ij} :=\left( X_i X_j u + X_j X_i u\right) /2 ,
\]
and $H$ is a given Hamiltonian, at least continuous and nondecreasing in $u$. Our main examples are the vector fields that generate the homogeneous Carnot groups \cite{Bell, BLU}, and in that case $D_{\mathcal X}u$  and $D^2_{\mathcal X}u$ are called, respectively, the horizontal gradient and  the horizontal Hessian.

A theory of fully nonlinear subelliptic equations was started recently by Bieske \cite{B1, B2} and Manfredi   \cite{Man, BBM}, and Monge-Amp\`ere equations of the form  \eqref{MAh} with $H=f(x)$ are listed among the main examples. For such equations on the Heisenberg group Gutierrez and Montanari \cite{GutMon} proved, among other things, a comparison principle among smooth sub- and supersolutions (see also \cite{GarTou} for related results).
% of the Monge-Amp\`ere equation  \eqref{MAh} with $H=f(x)$  .
%We are motivated by the recent interest of several researchers in developing a theory of subelliptic fully nonlinear equations.
%We are motivated by the recent interest for  fully nonlinear subelliptic equations.
%Moreover, a number of authors studied in the last five years several notions of convexity in Carnot groups \cite{LMS, DGN1, GutMon, BaRic, W2, Mag, Ric, GarTou, JLMS}, and one of their motivations is the connection with Monge-Amp\`ere equations. We will employ the notion of v-convexity introduced by Lu, Manfredi, and Stroffolini \cite{LMS} because it perfectly fits the viscosity methods we use throughout the paper, and it is the only one that makes sense for general vector fields. In the case of Carnot groups it coincides with the more geometric notion of {\em horizontal convexity}, a fact proved under different assumptions by several authors  \cite{LMS, BaRic, W2, Mag, Ric, JLMS}, see also \cite{DGN1}  for comparison with other notions. Here is the definition we use: an u.s.c. function $u$ on $\ol\Omega$ is convex with respect to the vector fields $X_1, ..., X_m$, briefly {\em ${\mathcal X}$-convex}, if it satisfies
%\[
%$
%-D^2_{\mathcal X}u\leq 0% \quad \mbox{in } \Omega
%$
%%\]
%in $\Omega$ in viscosity sense, that is, 
%\begin{equation}
%\label{hconv-test}
%D^2_{\mathcal X}\vf(x)\geq 0 \quad \forall\; \vf\in C^2(\Omega), x\in\argmax(u-\vf).
%\end{equation}
%The main example we have in mind to motivate 
An example that motivates the dependence on the gradient $D_{\mathcal X}u$ in $H$ is 
the prescribed horizontal Gauss curvature equation in Carnot groups, as defined by
%a subelliptic analogue of the prescribed Gauss curvature equation. In fact, 
Danielli, Garofalo and Nhieu \cite{DGN1},
%defined the horizontal Gauss curvature of the graph of the smooth function $u$ on a Carnot group as
%\[
%K_h(x) := \det (D_{\mathcal X}^2 u) \left(1+|D_{\mathcal X} u|^2\right)^{-\frac{m+2}{2}},
%\]
%so \eqref{MAh} becomes the prescribed horizontal Gauss curvature equation if 
%\begin{equation}
%\label{Gauss}
%H(x, r, q) = k(x) \left(1+|q|^2\right)^{\frac{m+2}{2}}
%\end{equation}
\begin{equation}
\label{Gauss:eq}
- \det (D_ {\mathcal X}^2 u) + k(x) \left(1+|D_{\mathcal X}u|^2\right)^{\frac{m+2}{2}}=0, \ \mbox{in}\ \Omega,
\end{equation}
for a given continuous $k : \ol\Omega \to ]0, +\infty[$.
%Finally, Gutierrez and Montanari \cite{GutMon} proved a comparison principle among smooth sub- and supersolutions of the Monge-Amp\`ere equation  \eqref{MAh} with $H=f(x)$  on the Heisenberg group.

In this paper we begin a study of the subelliptic Monge-Amp\`ere-type equations \eqref{MAh} within the theory of viscosity solutions. We present two comparison results that extend to the subelliptic setting a theorem of H. Ishii and P.-L. Lions for euclidean Monge-Amp\`ere equations \cite{IL}
(i.e., the case when the vector fields are the canonical basis of $\R^n$).
For the large literature on this case we refer to 
%we mean the case  It has a huge literature, especially on the regularity of solutions, see 
the recent surveys \cite{Ca, Tru06} and the references therein. The new difficulties we encounter are three.

1.  The PDE \eqref{MAh} is degenerate elliptic only on functions that are convex  with respect to the vector fields $X_1, ..., X_m$, briefly {\em ${\mathcal X}$-convex}. Following Lu, Manfredi, and Stroffolini \cite{LMS}  such a function is an u.s.c.  $u : \ol\Omega\to \R$ such that
 $
-D^2_{\mathcal X}u\leq 0% \quad \mbox{in } \Omega
$
in $\Omega$ in viscosity sense, that is, 
\begin{equation}
\label{hconv-test}
D^2_{\mathcal X}\vf(x)\geq 0 \quad \forall\; \vf\in C^2(\Omega), x\in\argmax(u-\vf).
\end{equation}
We refer to  the survey in \cite{BLU} for the %large 
recent literature on the notions of convexity in Carnot groups. %, in particular  the {\em horizontal convexity}.
%Moreover, a number of authors studied in the last five years several notions of convexity in Carnot groups \cite{LMS, DGN1, GutMon, BaRic, W2, Mag, Ric, GarTou, JLMS}, and one of their motivations is the connection with Monge-Amp\`ere equations. We will employ the notion of v-convexity introduced by because it perfectly fits the viscosity methods we use throughout the paper, and it is the only one that makes sense for general vector fields. In the case of Carnot groups it coincides with the more geometric notion of {\em horizontal convexity}, a fact proved under different assumptions by several authors  \cite{LMS, BaRic, W2, Mag, Ric, JLMS}, see also \cite{DGN1}  for comparison with other notions. 
Since ${\mathcal X}$-convex functions are not Lipschitz continuous, in general, 
%some results will restrict to 
we get better results in Carnot groups, where they are Lipschitz with respect to the intrinsic metric \cite{LMS, DGN1, Mag, Ric, JLMS}.

2. The operator in  \eqref{MAh} does not satisfy in general the standard structure conditions in viscosity theory. Therefore we consider equations of  the form
\begin{equation}
\label{MAlog}
-\log \det (D_{\mathcal X}^2 u)+ K(x, u, Du, D^2 u)=0, \ \mbox{in}\ \Omega,
\end{equation}
that verify the Lipschitz-type condition %(3.14)
with respect to $x$ of \cite{CIL} for {\em uniformly ${\mathcal X}$-convex} subsolutions.
%To overcome this problem we take the $\log$ of both terms in \eqref{MAh} and show that the new equation  i.e., functions such that, for a $\gamma>0$, 
%%\[
%$
%-D^2_{\mathcal X}u + \gamma I \leq 0
%% \quad \mbox{in } \Omega
%$
%%\]
%in $\Omega$ in viscosity sense. 
Our first main results states the comparison among semicontinuous sub- and supersolutions of this equation provided that either $K$ is strictly increasing in $u$ or that the subsolution is strict. Here $K$ is any degenerate elliptic operator satisfying the structure conditions of \cite{CIL}.

3. To cover the case of $H$ not strictly increasing in $u$, which is the most frequent in applications, we need to perturb a ${\mathcal X}$-convex subsolution to a uniformly ${\mathcal X}$-convex strict subsolution. 
%This was done in the euclidean case in \cite{IL} and we adapted the method to several subelliptic equations in  \cite{BM}. We are able to perform this construction under an additional condition on the 
%We adapt the method of \cite{IL} and \cite{BM} 
In the case of vector fields that generate a Carnot group we adapt the method of \cite{IL} and \cite{BM} 
 to get the following Comparison Principle, under essentially the same assumptions as the euclidean result %in
 of Ishii and Lions  \cite{IL}. 

\begin{theorem}
\label{teo:main}
Assume $H : \Omega\times\R\times \R^m \to ]0, +\infty[$ is continuous, nondecreasing in the second entry, and for all $R>0$ there is $L_R$ such that
\begin{equation}
\label{H1/m}
| H^{1/m}(x,r,q+q_1) -  H^{1/m}(x,r,q) | \leq L_R|q_1| %\mbox{ for all } 
\quad\forall \; x\in\ol\Omega, |r|\leq R, |q|\leq R, |q_1|\leq 1.
\end{equation}
Suppose the vector fields $X_1,...,X_m$ are the generators of a Carnot group on $\R^n$.
Let $u\,:\,\ol\Omega\to\R$ be a bounded, ${\mathcal X}$-convex, u.s.c. subsolution of \eqref{MAh} and $v\,:\,\ol\Omega\to\R$ be a bounded l.s.c. supersolution of \eqref{MAh}. Then
\begin{equation}
\label{comparison}
\sup_\Omega(u-v) \leq  \max_{\de\Omega}(u-v)^+.
\end{equation}
In particular, there is at most one  ${\mathcal X}$-%horizontally 
convex viscosity solution of \eqref{MAh} with prescribed continuous boundary data.
\end{theorem}
Note that it applies to the prescribed horizontal Gauss curvature equation %\eqref{MAh} 
\eqref{Gauss:eq}.

In Section \ref{uno} we state the Comparison Principle for the equation \eqref{MAlog} 
% and outline 
with the main lemma needed for its proof. Section \ref{due} is devoted to recalling the definition of generators of a Carnot group and stating a few facts about them. Finally, in Section \ref{tre} we outline the construction of the strict subsolution and the rest of the proof of Theorem \ref{teo:main}. Our paper \cite{BM3} contains the full proofs of these results, some extensions and variants, the existence of solutions to the Dirichlet problem via the Perron-Ishii method, and further examples and bibliography.

%Other subelliptic fully nonlinear equations were studied by Bieske \cite{B1, B2} and ourselves \cite{BM}, see also the references therein. 
%%

%%%%%%%%%%%%%%%%%%%%%%%%%%%%%%%%%%%%%%%%%%%%%%%%%%%
%%    Section 1
%%%%%%%%%%%%%%%%%%%
\section{Definitions and comparison with strict subsolutions}
\label{uno}
Let $\sigma$ be the $n\times m$ matrix-valued function whose columns  $\sigma^j$ are the coefficients of the vector fields $X_1,...,X_m$, $j=1,\cdots, m$. %Throughout the paper w
We assume $\sigma^j_i=\sigma_{ij} \in C^{1,1}(\ol\Omega)$ for all $i, j$. Observe that, for a smooth function $u$,
$D_{\mathcal X}u(x) = \s(x)^TDu(x)$ and
\[
%\begin{equation}
%\label{D2h}
D_{\mathcal X}^2 u(x) = \sigma^T(x) D^2 u(x)\, \sigma(x) + Q(x,Du), \; Q_{ij}(x, p):= \left[D\sigma^j\,\sigma^i+D\sigma^i\,\sigma^j\right](x)\cdot \frac{p}{2}.
%\end{equation}
\]
Therefore we rewrite \eqref{MAh} and \eqref{MAlog} in the form $G(x,u, Du, D^2u)=0$ with $G$ proper
in the sense of \cite{CIL}.

%the usual definition of viscosity subsolution \cite{CIL} makes sense for equations of the form \eqref{MAh} and \eqref{MAlog}. 
%The definition of supersolution $v$ has to be modified as in \cite{IL}, Section V.3: we require that
%\[
%%$
%-\log\det (D^2_{\mathcal X}\vf)
%%\sigma^T(x) D^2 u\, \sigma(x)+Q(x, Du))
%+ K(x, v, D\vf, D^2\vf)\geq 0\quad\forall x\in \arg\min (v-\vf),
%%$  
%\]
%%at all minimum points $x$ of $v-\vf$,
%for any $\vf\in C^2(\Omega)$ such that $D^2_{\mathcal X}\vf(x)$ is positive definite.
%%where $\vf$ is any function in $C^2(\Omega)$ satisfying in addition  $D^2_{\mathcal X}\vf(x)\geq 0$.

We say that a continuous function $F : \ol \Omega\times \re\times \re^n\times S^n   \rightarrow \re$ {\em satisfies the structure conditions} (of viscosity theory) on a given set of admissible symmetric matrices ${\mathcal M} \subseteq S^n$ if it is nondecreasing in the second entry, nonincreasing in the last entry for matrices in ${\mathcal M}$, and for some modulus $\omega$
\[
F\left(y, r, \frac {x-y}{\epsilon},Y\right) - F\left(x, r,  \frac {x-y}{\epsilon},X\right)\leq \omega\left( |x-y|\left(1+\frac {|x-y|}{\epsilon}\right)\right)
\]
for all $\epsilon>0$, $x, y\in \ol\Omega$, $r\in \R$, $X,Y\in {\mathcal M}$ satisfying  
\begin{eqnarray}\label{STRUT}
&& \ -\frac{3}{\epsilon}\left(\begin{array}{cc}
I&0\\
0&I\end{array}\right)\leq
\left(\begin{array}{cc}
X&0\\
0&-Y\end{array}\right)\leq
\frac{3}{\epsilon}\left(\begin{array}{cc}
I&-I\\
-I&I\end{array}\right).\label{structure}
\end{eqnarray}

%\begin{definition}
We say that $u:\ol\Omega\to\R$ u.s.c.  is
%$ \rm{convex in $\Omega$ with respect to the fields $X_1,\cdots,X_m$}, or \rm{horizontally convex},
{\em uniformly ${\mathcal X}$-convex} if for some $\gamma>0$
\begin{equation}\label{unifconvex}
D_{\mathcal X}^2 \vf \geq \gamma I, \quad \forall\; \vf\in C^2(\Omega), x\in\argmax(u-\vf), %\ in\ \Omega,
\end{equation}
where $I$ denotes the identity matrix. In other words, with the notations of \cite{CIL}, $(p,X)\in \mathcal J^{2,+}u(x)$ satisfies
%\[
$
 \sigma^T(x) X\, \sigma(x) + Q(x,p)  \geq \gamma I,
$
%\]
and this inequality defines the set of admissible matrices ${\mathcal M}={\mathcal M}(p,\gamma)$.

The main ingredient for the results
%Comparison Principle
 of this section is the following.
\begin{lemma}\label{lemmadiff}
For each $\gamma>0$ the function
$
F(x,p,X):= -\log \det ( \sigma^T(x) X\, \sigma(x)+Q(x,p))
$
satisfies the structure conditions on ${\mathcal M}(p,\gamma)$.
%Then for all $\gamma>0$ there is a constant $C>0$ such that
%\begin{eqnarray}
%&&F(y,\frac{x-y}{\epsilon}, Y)-F(x,\frac{x-y}{\epsilon}, X)\leq C(|x-y|+\frac{|x-y|^2}{\epsilon})\label{diff}
%\end{eqnarray}
% for all  $X, Y \in S^n$ satisfying \ref{structure} and  $p\in \R^n$ such that
%\begin{equation}
% \sigma^T(x) X\, \sigma(x)+Q(x,p)\geq \gamma I,\  \sigma^T(y) Y\, \sigma(y)+Q(y,p)\geq \gamma I.\label{gamma}
%\end{equation}
\end{lemma}
%\begin{proof} The degenerate ellipticity is standard. Next we  represent  $F$  as a maximum of  operators that satisfy the structure conditions. This can be done by the following formula, holding for $A\in S^m$, $A\geq \gamma\, I$,
%%\begin{eqnarray}&&
%\[
%\log\text{det}(A)=%\label{represlog}\\
%%&&=
%\min\{ m\log a-m+\mbox{tr}(AM) : a>0, M\in S^m, 0\leq M\leq \frac{1}{\gamma}I, \text{det}\ M=a^{-m}\}.
%\]
%%\nn \end{eqnarray} 
%The rest of the proof follows the methods of \cite{CIL}, see \cite{BM3}.
%%Hamilton-Jacobi-Bellman
%\end{proof}
The proof relies on a  representation of  $F$  as a maximum of  operators that satisfy the structure conditions, via the following formula, holding for $A\in S^m$, $A\geq \gamma\, I$,
%\begin{eqnarray}&&
\[
\log\text{det}(A)=%\label{represlog}\\
%&&=
\min\{ m\log a-m+\mbox{tr}(AM): a>0, M\in S^m, 0\leq M\leq \frac{1}{\gamma}I, \text{det}\ M=a^{-m}\}.
\]

The next two Comparison Principles can now be proved by standard methods in viscosity theory  \cite{CIL}, see \cite{BM3} for the details.
%The definition of supersolution $v$ has to be modified as in \cite{IL}, Section V.3: 
In the definition of {\em supersolution} $v$ of \eqref{MAh} and \eqref{MAlog} we restrict to ${\mathcal X}$-convex test functions. %, i.e.,
E.g., for \eqref{MAlog} we require that
\[
%$
-\log\det (D^2_{\mathcal X}\vf)
%\sigma^T(x) D^2 u\, \sigma(x)+Q(x, Du))
+ K(x, v, D\vf, D^2\vf)\geq 0\quad\forall x\in \arg\min (v-\vf),
%$  
\]
%at all minimum points $x$ of $v-\vf$,
for any $\vf\in C^2(\Omega)$ such that $D^2_{\mathcal X}\vf(x)$ is positive definite, cf. \cite{IL}, Section V.3.
%where $\vf$ is any function in $C^2(\Omega)$ satisfying in addition  $D^2_{\mathcal X}\vf(x)\geq 0$.

 \begin{theorem}
\label{confronto} 
Assume $K : \ol \Omega\times \re\times \re^n\times S^n   \rightarrow \re$ satisfies the structure conditions on $S^n$.
%(\ref{B}),  (\ref{STRUT}), (\ref{sigma}),  (\ref{Q0}).
Let $u\,:\,\ol\Omega\to\R$ be bounded, uniformly ${\mathcal X}$-convex, and  for all open $\Omega_1$ with $\ol \Omega_1\subseteq \Omega$ there is $\gamma_1>0$ such that $u$ is subsolution of
\begin{equation}
\label{MAlogstretta}
-\log\det (D^2_{\mathcal X}u)
%\sigma^T(x) D^2 u\, \sigma(x)+Q(x, Du))
+ K(x, u, Du, D^2u)\leq -\gamma_1,\ \text{in}\ \Omega_1.
\end{equation}
Let $v\,:\,\ol\Omega\to\R$ be a bounded l.s.c. supersolution of \eqref{MAlog}. Then
%\begin{equation}
%\label{comparison}
$$\sup_{\Omega}(u-v)\leq \max_{\partial\Omega}(u-v)^+ .$$
\end{theorem}
\begin{theorem}
\label{confronto2} 
The conclusion of the previous theorem remains true if $u$ is a  subsolution of \eqref{MAlog}, not necessarily strict, provided that, for some $C>0$,
$$K(x,r,p,X)-K(x,s,p,X)\geq C(r-s), \; -M\leq s\leq r\leq M,\; M:=\max %\{\sup|u|, \sup|v|\}.$$
\{\|u\|_\infty, \|v\|_\infty\}.$$
Under this condition there is at most one uniformly ${\mathcal X}$-convex viscosity solution of \eqref{MAlog} with prescribed continuous boundary data.
\end{theorem}

%(in viscosity sense), where
%$
%Q_{ij}(x, p)= \left(D\sigma^j(x)\,\sigma^i(x)+D\sigma^i(x)\,\sigma^j(x)\right)\cdot \frac{p}{2},
%$
%\end{definition}
%
%Note that, for smooth $u$, the inequality can be written as
%$$D_{\mathcal H}2 u\geq 0,\ in \Omega,$$
%where
%$D_{\mathcal H}2 u$ is the symmetrized Hessian matrix with respect to the fileds, i.e.,
%$(D_{\mathcal H}2 u)_{ij} = \frac{X_i(X_ju)+ X_j(X_iu)}{2}$.
%Therefore, if $X_1,\cdots,X_m$ are the generators of a Carnot group $\mathcal{G}$, our definition coincides with the definition of
%convexity in $\mathcal{G}$ in viscosity sense (v-convexity) of Lu, Manfredi, Stroffolini (\cite{LMS}).
%A more geometric notion of convexity in $\mathcal{G}$, called horizontal convexity, was introduced and studied
%in the same seminal paper (\cite{LMS}) and, independently, by Danielli, Garofalo, and Nhieu  (\cite{DGN}).
%The equivalence of the two notions was proved by several authors, first in the Heisenberg groups (\cite{LMS}),
%\cite{BaRi}, and then in general Carnot groups [Wang, Magnani, JLMS].
%This justifies the terminology of our definition, although we also consider vector fields that are not associated with a group structure.

%%%%%%%%%%%%%%%%%%%%%%%%%%%%%%%%%%%%%%%%%%%%%%%%%%%
%%    Section 2
%%%%%%%%%%%%%%%%%%%

\section{Generators of Carnot groups}
\label{due}
We begin with recalling some well-known definitions. We adopt the terminology and notations of the recent book \cite{BLU}.
Consider a group  operation $\circ$ on $\R^n=\R^{n_1}\times...\times\R^{n_r}$ with identity $0$, such that $(x,y)\mapsto y^{-1}\circ x$ is smooth,  and the dilation $\d_\l\,:\,\R^n\to\R^n$ 
\[
\d_\l(x) = \d_\l(x^{(1)},...,x^{(r)}) := (\l x^{(1)}, \l^2 x^{(2)},..., \l^r x^{(r)}), \quad x^{(i)}\in \R^{n_i}.
\]
If $\d_\l$ is an automorphism of the group $(\R^n, \circ)$ for all $\l>0$, %then 
$(\R^n, \circ, \d_\l)$ is a homogeneous Lie group on $\R^n$. We say that $m=n_1$ smooth vector fields $X_1, ..., X_m$ on $\R^n$ generate $(\R^n, \circ, \d_\l)$, and that this is a (homogeneous) Carnot group, if
%\begin{itemize}
%\item 
$X_1, ..., X_m$ are invariant with respect to  the left translations on $\R^n$ $\tau_\a(x):= \a\circ x$ for all $\a\in\R^n$,
%\item 
$X_i(0) = \de /\de x_i$, $i=1,...,m$,
%\item 
and the  rank  of the Lie algebra generated by $X_1, ..., X_m$ is $n$ at every point $x\in\R^n$.
%\end{itemize}
%
We refer, e.g., to \cite{Bell, BLU} for the connections of this  definition with the classical one in the context of  abstract Lie groups and for the properties of the generators. We will use only  the following property, and refer to Remark 1.4.6, p. 59 of \cite{BLU} for more precise informations.

\begin{proposition}
\label{Carnot}
 If  $X_1, ..., X_m$ are generators of a Carnot group, then
\[
X_j(x) = \frac{\de}{\de x_j} + \sum_{i=m+1}^n \s_{ij}(x) \frac{\de}{\de x_i}
\]
with $\s_{ij}(x)=\s_{ij}(x_1,...,x_{i-1})$ homogeneous polynomials of a degree $\leq n-m$.
\end{proposition}
For generators of Carnot groups %various authors proved %, under different assumptions, 
the Lipschitz
continuity of ${\mathcal X}$-convex functions with respect to the intrinsic metric %of the group 
and bounds on the horizontal gradient in the sense of distributions were studied in \cite{LMS, DGN1, Mag, Ric, JLMS}.
%From those results we obtain 
We deduce the following gradient bound in viscosity sense.
%%%%
% Da mettere capitolo Carnot
%%%%%%
\begin{proposition}\label{stimagradiente}
Let $u$ be convex in $\Omega$ with respect to the generators of a Carnot group. Then, for every open $\Omega_1$
with $\ol\Omega_1\subseteq \Omega$, there exists a constant $C$ such that $u$ is a viscosity subsolution of 
$ |\sigma^T(x)\, Du|\leq C$ %\  \mbox{in } \Omega_1.
in $\Omega_1$.
%in viscosity sense.
\end{proposition}
%\begin{proof}
%Since $u$ is u.s.c. it is locally bounded from above. Then $h$-convexity implies local Lipschitz continuity with respect to Carnot-Caratheodory distance, by a result of Magnani [ ] and Rickly [  ], see also [JLMS].......\end{proof} 

%%%%%%%%%%%%%%%%%%%%%%%%%%%%%%%%%%%%%%%%%%%%%%%%%%%
%%    Section 3
%%%%%%%%%%%%%%%%%%%

\section{Outline of proof of the Comparison Principle}
\label{tre}
In this section we outline the proof of Theorem \ref{teo:main}. 
The definition of viscosity supersolution of \eqref{MAh} uses only ${\mathcal X}$-convex test functions as  in Section \ref{uno} and in \cite{IL}.
Given the ${\mathcal X}$-convex subsolution $u$ we consider
$$
u_{\epsilon,\mu}(x):=u(x)+\epsilon e^{\mu\frac{\sum_{i=1}^m |x_i|^2}{2}},
$$
for positive $\epsilon, \mu$. %, \lambda$.
%for small $\epsilon>0$. For $\mu>0$ it is easy to see  and a suitable choice of the other parameters. 
A calculation using Proposition \ref{Carnot} shows that $u_{\epsilon,\mu}$ is uniformly ${\mathcal X}$-convex with $\gamma = \epsilon\mu$.
\begin{lemma}
For any open $\Omega_1$ with $\ol \Omega_1\subseteq \Omega$ there are positive constants 
$\ol \mu$, independent of $\epsilon$, and $\gamma_2$ 
such that, for $\mu\geq\ol\mu$, $u_{\epsilon,\mu}$ is a subsolution of
\begin{equation*}
-{\rm det}^{1/m} (D^2_{\mathcal X}u) + 
H^{1/m}(x, u, D_{\mathcal X}u)\leq -\gamma_2,\ \text{in}\ \Omega_1.
\end{equation*}
\end{lemma}
The proof of the lemma relies on the Minkowski inequality for %(\cite{HJ})
$\text{det}^{1/m}(A+B)$ %with $A>0$, $B \geq 0$, 
with $A, B$ positive definite, and the identity
%\[
%\begin{equation}\label{NoDa}
$
\det(I+q q^T)= 1+|q|^2
$ %\ q \mbox{ is a } m\times 1 \mbox{  column vector},
%\end{equation}\]
%where $(q\otimes q )_{ij}:=q_iq_j$. 
for any column vector $q\in\R^m$. Moreover, by the boundedness of $u$ and $D_{\mathcal X}u$, Proposition \ref{stimagradiente},  we can assume the Lipschitz property \eqref{H1/m} with a uniform constant $L_R$. The rest of the proof goes along the lines of \cite{BM}.

Next, we exploit again the boundedness of  $D_{\mathcal X}u$ in viscosity sense to see that 
$u_{\epsilon,\mu}$ satisfies \eqref{MAlogstretta} for $K=\log H$ and a suitable $\gamma_1>0$. Then Theorem \ref{confronto} applies and gives $\sup_{\Omega}(u_{\epsilon,\mu}-v)\leq \max_{\partial\Omega}(u_{\epsilon,\mu}-v)^+ .$ Letting ${\epsilon}\to 0$ gives the conclusion.

%\begin{remark}{\rm
%Our definition of solution to \eqref{MAh}  is equivalent to the one proposed in \cite{DGN1}, but the notions of sub- and supersolution are weaker because we do not assume the continuity. This is useful in the construction of solutions to the Dirichlet problem and in problems of approximation, homogenization, etc., where the relaxed viscosity semilimits are used.
%}
%\end{remark}

%\smallskip

%%%%%%%%%%%%%%%%%%%%%%%%%%%%%%%%%%%%
\medskip
\noindent
{\sc {\bf Acknowledgments.}} 
%\noindent {\sc Acknowledgments } 
The authors are grateful to Roberto Monti and Luigi Salce for several useful talks.
%the proof of Lemma (\ref{salce}).
Work partially supported by 
the Italian M.I.U.R. project "Viscosity,
metric, and control theoretic methods for nonlinear partial
differential equations''.

%\newpage

%%%%%%%%%%%%%%%%%%%%%%%%%%%%%%%%%%%%%%%%%%%%%%%%%%%\

%\nocite{Fri:64}
%\bibliography{fp}
%\bibliographystyle{plain}

%\vskip 0.3truecm

%\noindent {\bf Address of the authors}\\
%Dipartimento di Matematica Pura e Applicata,\\
%Universit\`a degli Studi di Padova,\\
%Via Trieste, 63, 35121, Padova, Italy,\\
%bardi@math.unipd.it, mannucci@math.unipd.it\\

\end{document}